\documentclass[letterpaper, 10 pt, conference]{ieeeconf}
\usepackage{amssymb}
\usepackage{color}
\usepackage{bbm}

\usepackage{graphicx}
\usepackage{mathtools}
\usepackage{multirow}
\usepackage{lineno,hyperref}
\usepackage{url}
\usepackage{cleveref}
\usepackage{subcaption}
\usepackage{bm}
\usepackage{tikz}
\usetikzlibrary{arrows.meta}
\usepackage{algorithm}
\usepackage{algorithmic}
\usepackage{caption}
\usepackage{makecell}

\usepackage{mathtools}
\mathtoolsset{showonlyrefs,showmanualtags}

\DeclareMathAlphabet{\pazocal}{OMS}{zplm}{m}{n}

\newcommand{\ve}{\boldsymbol}

\renewcommand{\geq}{\geqslant}
\renewcommand{\epsilon}{\varepsilon}

\usepackage{pgfplots}
\DeclareUnicodeCharacter{2212}{−}
\usepgfplotslibrary{groupplots,dateplot}
\usetikzlibrary{patterns,shapes.arrows}

\modulolinenumbers[5]

\newtheorem{remark}{Remark}[section]

\usepackage{multirow}
\usepackage{color}
\usepackage{paralist}
\usepackage[space]{cite}

\pgfplotsset{compat=1.18}

\IEEEoverridecommandlockouts                        

\title{\LARGE \bf
UQ-ACE: Stability-Aware Counterfactual Controller Tuning via Bayesian SINDy
}

\author{Margarita A. Guerrero$^{1}$, Rodrigo A. González$^{2}$,
and Cristian R. Rojas$^{1}$%
\thanks{This work has been partially supported by the Swedish Research Council under contract number 2023-05170, and by the Wallenberg AI, Autonomous Systems and Software Program (WASP) funded by the Knut and Alice Wallenberg Foundation.}%
\thanks{$^{1}$Margarita A. Guerrero and Cristian R. Rojas are with the
Department of Decision and Control Systems, KTH Royal Institute of
Technology, 100 44 Stockholm, Sweden
(e-mails: mags3@kth.se, crro@kth.se).}%
\thanks{$^{2}$Rodrigo A. González is with the Department of Mechanical
Engineering, Eindhoven University of Technology,
Eindhoven, The Netherlands
(e-mail: r.a.gonzalez@tue.nl).}%
}

\begin{document}

\maketitle

\begin{abstract}
Retuning feedback controllers in critical processes is challenging when experiments are costly, open-loop excitation is infeasible, and candidate controllers must be screened before deployment. The existing Adaptive Sampling for Counterfactual Explanations (ACE) method provides a sample-efficient approach for finding minimal controller-parameter changes that recover prescribed closed-loop performance, but it does not explicitly account for stability of candidate controllers proposed during the search, which may be destabilizing under plant uncertainty. This paper introduces UQ-ACE, a stability-aware retuning method based on Bayesian SINDy. In UQ-ACE, closed-loop trajectories are used to infer a posterior distribution over local dynamical models, from which candidate controllers are screened for stability before implementation. Extensive simulations show that UQ-ACE achieves valid retunings with few closed-loop experiments while preventing deployment of unstable candidates.
\end{abstract}

\begin{keywords}
Bayesian Methods; Learning for Control; Iterative Control Tuning. 
\end{keywords}

\section{Introduction}

In setpoint-regulation problems, feedback controllers maintain process variables at desired values despite disturbances and changes in operating conditions. In large-scale and safety-critical industrial systems, retuning cannot rely on trial-and-error: opening the loop or applying aggressive excitation is often infeasible because interruptions are costly and unsafe operating regimes may be reached \cite{aastrom_05,khosravi_23}. At the same time, industrial plants often store long-horizon closed-loop records in historical servers, including process
variables, setpoints, manipulated variables, and controller modes
\cite{crompton2021}. These records provide prior information
about feasible operating regions, but numerical controller parameters, such as
PID gains or lead--lag time constants, are often not logged and may need to be
inferred from data before retuning. This motivates data-driven retuning methods
that exploit historical closed-loop data while requiring only a small number of
additional safe experiments.

The objective in controller retuning is not only to improve closed-loop
performance, but also to do so safely. Here, safe tuning refers to a
data-driven retuning procedure that avoids controller-parameter updates likely
to produce unstable or otherwise unsafe closed-loop behavior, such as excessive
oscillations, actuator saturation, large overshoot, or loss of stability. This view is
aligned with safe learning in control, where data-driven performance
improvement must be balanced against constraints, uncertainty, and safety
requirements~\cite{hewing_20}. Several safe learning and data-driven optimization approaches address related
issues: SafeOpt~\cite{sui2015safe} and safe Bayesian optimization methods restrict the
search to parameters that satisfy a safety condition with high probability
\cite{berkenkamp2016safe,baumann2025safety}, while constrained Bayesian
optimization has been used for safety-aware cascade controller tuning from
closed-loop experiments~\cite{khosravi_23}. Similar ideas also appear in
learning-based MPC and safe learning in control, where performance is improved
from data while accounting for constraints, uncertainty, and safety
requirements~\cite{hewing_20,Wabersich_23}. However, these methods typically encode safety
through performance thresholds, constraint surrogates, safe sets, or
model-based uncertainty descriptions. In contrast, our goal is to retune an
existing controller by finding a minimal and actionable parameter change.

This minimal-change retuning objective naturally connects to counterfactual explanations (CFE)~\cite{wachter2018counterfactual}, which seek the smallest actionable change that flips a predicted outcome. Adaptive Sampling for Counterfactual Explanations (ACE)~\cite{guerrero_26} provides a promising approach to this problem. In the
controller-retuning setting, ACE casts the performance requirement as a
pass/fail classification task over the controller-parameter space. Starting
from a failing baseline controller, ACE uses a Gaussian-process classifier (GPC) and
an acquisition function to search for a nearby controller that flips the label
from fail to pass with as few closed-loop experiments as possible~\cite{guerrero_26_ifac}.
Despite its sample efficiency, ACE does not explicitly quantify how uncertainty in the local closed-loop dynamics affects the stability of a candidate controller before deployment. Consequently, a candidate that appears promising from the pass/fail classifier may still lead to undesirable or unstable closed-loop behavior when tested.

Motivated by this limitation, we propose UQ-ACE, a stability-aware data-driven
controller tuning method. UQ-ACE extends ACE by adding a stability screen based
on Bayesian sparse
identication of nonlinear dynamics (SINDy)~\cite{brunton_16,hirsh_22}. The closed-loop trajectories collected
during tuning are used to infer a posterior distribution over local dynamical
models. Each posterior sample induces a possible local closed-loop model for a
candidate controller, and the candidate is accepted for testing only if its
posterior stability probability exceeds a prescribed threshold. In this way, UQ-ACE preserves the
sample-efficient and minimal-change nature of ACE while adding an
uncertainty-aware stability requirement before deployment.

The main contributions of this paper are as follows:
\begin{itemize}
    \item We formulate stability-aware counterfactual controller tuning as a
    data-driven retuning problem with performance, minimal-change, and posterior
    stability requirements.
    \item We introduce UQ-ACE, which augments counterfactual controller tuning with a Bayesian SINDy stability screen
    built from closed-loop tuning data.

    \item We evaluate the proposed method on an axial
    compressor--plenum--throttle system and compare it with the original ACE
    procedure without posterior stability screening.
\end{itemize}
The remainder of the paper is organized as follows.
Section~\ref{sec:problem} formulates the safe retuning problem, and
Section~\ref{sec:preliminaries} reviews ACE and introduces the Bayesian SINDy
model. Section~\ref{sec:approach} presents the proposed UQ-ACE method, while
Section~\ref{sec:simulations} reports the simulation results.
Section~\ref{sec:conclusions} concludes the paper.

\section{Setup and Problem Formulation}
\label{sec:problem}
We consider the retuning of critical feedback loops for which open-loop experiments are infeasible and aggressive excitation is not acceptable. The goal is to derive an iterative
procedure that produces the smallest possible change in the baseline controller parameters while maintaining safe operation and ensuring control performance.

Consider the discrete-time closed-loop system described in state-space form by
\begin{subequations}\label{eq:System}
\begin{align}
    \ve x_{t+1} &= \ve f(\ve x_t,\ve u_t)+\ve w_t, \label{eq:System_a}\\
    \ve x_{c,t+1}&=\ve A_c(\ve\theta)\ve x_{c,t}+\ve B_c(\ve\theta) \ve e_t, \label{eq:System_b} \\
    \ve u_t&=\ve C_c(\ve\theta)\ve x_{c,t}+\ve D_c(\ve\theta) \ve e_t , \label{eq:System_c} \\
    \ve e_t&=\ve r_t-\ve C_p \ve x_t, \label{eq:System_d}
\end{align}
\end{subequations}
where $\ve x_t\in \mathbb{R}^{n_x},\ve u_t \in \mathbb{R}^{n_u}, \ve r_t \in \mathbb{R}^{n_r}$ and $\ve e_t \in \mathbb{R}^{n_r}$ are the state, input, reference and reference error signals respectively, and $\ve w_t$ is the process noise, which is assumed white and Gaussian, $\ve w_t \sim \mathcal{N}(\mathbf{0},\bm{\Sigma}_w)$. The transition map $\ve f$ is unknown and possibly nonlinear, and $(\ve{A}_c(\ve\theta),\ve{B}_c(\ve\theta),\ve{C}_c(\ve\theta),\ve{D}_c(\ve\theta))$ is the linear and time-invariant state-space representation of the controller, whose parameters are contained in a vector $\ve \theta \in \Theta$. The state $\ve x_t$ is assumed to be available, and the matrix $\ve C_p$ selects the measured components of the state for control.

Let $\ve\theta_0$ denote the baseline controller currently operating in closed loop. 
We assume the availability of batches of closed-loop data collected under $m$ 
previously tested controller parameters $\ve\theta_1,\ldots,\ve\theta_m$. For each 
controller $\ve\theta_i$, the recorded data are given by
\begin{equation}
    \mathcal{D}_i
    =
    \{\ve r_t^{\langle i\rangle },\ve u_{t}^{\langle i\rangle },\ve x_{t}^{\langle i\rangle }\}_{t=0}^{N_i},
    \qquad i=1,\ldots,m,
    \label{eq:experiments}
\end{equation}
where $N_i$ denotes the length of the $i$th experiment. These trajectories are recorded under controllers operated within an admissible closed-loop envelope; they therefore provide prior information about feasible controller parameters and locally safe dynamic behavior.
\begin{remark}
    If the parameters $\ve\theta_i$ are not logged, they can be recovered exactly by fitting \eqref{eq:System_b}-\eqref{eq:System_d} to $\mathcal{D}_i$, provided the controller parametrization is known and the data are sufficiently rich for $\ve\theta_i$ to be uniquely determined. If the input measurements are noisy, $\ve\theta_i$ can instead be estimated by nonlinear least squares or subspace methods \cite{Ljung1999SystemUser}. These estimates can then be used in the proposed retuning procedure.
\end{remark}
We seek to construct, through a sequence of safe closed-loop trials informed by the previously collected data, a controller parameter vector $\ve\theta$ that remains as close as possible to $\ve\theta_0$ under a prescribed metric $d$ on $\Theta$, while satisfying time-domain performance specifications. In addition, because the controller is to be deployed on a critical feedback loop, each candidate parameter vector must satisfy a stability-screening requirement under the uncertainty induced by the available closed-loop data. This leads to the safe retuning problem
\begin{equation}
\begin{aligned}
    \min_{\ve\theta\in\Theta}
    \; \; d(\ve\theta,\ve\theta_0) \; \;
    \mathrm{s.t.}
    \;\; & \textnormal{performance specifications},\\[-1.3ex]
    & \textnormal{closed-loop stability requirement}.
    \label{eq:optim_problem}
\end{aligned}
\end{equation}

The formulation in \eqref{eq:optim_problem} captures the desired retuning objective, but the performance and stability constraints are not available in closed form. We therefore address it via a sample-efficient iterative procedure that sequentially proposes candidate controllers near the operating condition and screens them using a data-driven stability test.

\section{Preliminaries}
\label{sec:preliminaries}

In this section, we introduce the technical preliminaries required to address the retuning problem \eqref{eq:optim_problem} via a counterfactual explanability framework.

\subsection{\hspace{-0.15cm}Adaptive\hspace{-0.01cm} sampling\hspace{-0.01cm} for\hspace{-0.01cm} counterfactual\hspace{-0.01cm} explanations\hspace{-0.04cm} (ACE)}
\label{sec:ace}

Consider a controller parameterized by $\bm{\theta}\in\Theta$, and let $h\colon\Theta\to\{0,1\}$ indicate whether a closed-loop experiment with $\bm{\theta}$ meets the prescribed performance specifications, corresponding to the first constraint in~\eqref{eq:optim_problem}. Specifically, $h(\bm{\theta})=1$ if the specifications are satisfied (pass), and $h(\bm{\theta})=0$ otherwise (fail). Given a baseline controller $\bm{\theta}_0$ with $h(\bm{\theta}_0)=0$, the minimal retuning meeting the performance specifications is a \emph{counterfactual explanation} (CFE)
\cite{wachter2018counterfactual} of $\bm{\theta}_0$:
\begin{equation}\label{eq:cfe}
  \min_{\bm{\theta}\in\Theta}\; d(\bm{\theta},\bm{\theta}_0)
  \quad \text{s.t.} \quad h(\bm{\theta})=1.
\end{equation}
Since the pass/fail map $h$ is only observed through costly closed-loop
experiments, ACE \cite{guerrero_26}, as adapted to controller retuning in \cite{guerrero_26_ifac}, solves \eqref{eq:cfe} with as few queries as possible.

\emph{Penalized formulation.}
Following \cite{rasmussen2006gaussian}, we model the pass/fail outcome as $h(\bm{\theta}) = \mathbb{H}\big(f_h(\bm{\theta}) - 0.5\big)$, where $\mathbb{H}$ denotes the Heaviside step function and $f_h$ is a smooth latent score. The decision boundary is then the level set $\{\bm{\theta}\colon f_h(\bm{\theta}) = 0.5\}$. Given a dataset of $m$ controller-label pairs $\mathcal{P} = \{(\bm{\theta}_i, h(\bm{\theta}_i))\}_{i=1}^{m}$, each containing a parameter vector with its observed pass/fail outcome, we replace
\eqref{eq:cfe} by
\begin{align}
  \min_{\bm{\theta}\in\mathcal{A}}\; J_k(\bm{\theta}) := {}&
    d(\bm{\theta},\bm{\theta}_0) + \lambda_k\,\big|f_h(\bm{\theta})-0.5\big| \nonumber\\
  &+ \beta\, g(\bm{\theta}-\bm{\theta}_0) + l(\bm{\theta};\mathcal{P}).
  \label{eq:ace_cost}
\end{align}
The first two terms keep $\bm{\theta}$ close to $\bm{\theta}_0$ and
drive it to the boundary, where the label flips with minimal parameter
change; $\lambda_k$ grows with the iteration index $k$, as in the
penalty method \cite{luenberger2008linear}. The set
$\mathcal{A}\subseteq\Theta$ fixes parameters that cannot be modified
(\emph{actionability}), $g=\|\cdot\|_1$ with $\beta>0$ favors retuning
few parameters (\emph{sparsity}), and
$l(\bm{\theta};\mathcal{P})\in\{0,+\infty\}$ discards controllers that
the local outlier factor \cite{breunig2000lof} flags as outliers of $\mathcal{P}$, keeping experiments within the observed operating
envelope (\emph{plausibility}).

\emph{Surrogate model.} As only labels are observed, $f_h$ is replaced by a Gaussian process classifier surrogate $\hat{f}_h$ fitted to $\mathcal{P}$. Under the Laplace approximation, the latent posterior at $\bm{\theta}$ is $\mathcal{N}(\mu_a(\bm{\theta}),\sigma_a^2(\bm{\theta}))$~\cite[Eq. (6)]{guerrero_26}, and the pass probability is approximated by \cite{rasmussen2006gaussian}
\begin{equation}
   \hat{p}_{h}(\bm{\theta})= \sigma\big(\mu_a(\bm{\theta})\,
  (1+\pi\sigma_a^2(\bm{\theta})/8)^{-1/2}\big)\in[0,1], \notag 
\end{equation}
where $\sigma$ denotes the logistic sigmoid.

\emph{Adaptive sampling.} The next controller to test maximizes the
expected improvement (EI) \cite{jones1998efficient}
\begin{equation}\label{eq:ei}
  \bm{\theta}_{m+1} = \underset{\bm{\theta}\in\mathcal{A}}{\arg \max}\;
  \mathbb{E}\big[\max\{J_k(\bm{\theta}^\ast)-J_k(\bm{\theta}),\,0\}\big],
\end{equation}
where $\bm{\theta}^\ast$ is the lowest-cost controller in $\mathcal{P}$ and the expectation is over the posterior of $\hat{f}_h$. The EI is computed by a Monte Carlo approximation with correlated samples of $(\hat{f}_h(\bm{\theta}),\hat{f}_h(\bm{\theta}^\ast))$ and maximized by L-BFGS-B \cite{byrd1995limited}; it balances exploitation (low predicted cost) and exploration (high uncertainty). After each experiment, $\mathcal{P}$ and $\hat{f}_h$ are updated and
$\lambda_{k+1}=\lambda_k^{p}$, $p>1$, until convergence and $\lambda_k\geq\lambda_{\max}$. The retuned controller is the point of a dense Sobol' set closest to $\bm{\theta}_0$ with $\hat{p}_{h}(\bm{\theta})\approx 0.5$; the procedure restarts until the controller passes and no closer one is found.

\begin{remark}\label{rem:dataset}
    ACE starts with $\mathcal{D}_i, \; i=1,\dots,m_{hist}$ (cf.~\eqref{eq:experiments}) historical data, from which each controller $\ve \theta_i$ is identified and its performance is evaluated to obtain $\mathcal{P}_{hist} = \{(\bm{\theta}_i, h(\bm{\theta}_i))\}_{i=1}^{m_{hist}}$. Historical controller-label pairs plus newly collected labeled controller candidates are defined as $\mathcal{P}=\mathcal{P}_{hist} \cup \mathcal{P}_{online}=\{(\ve{\theta}_i, h(\ve{\theta}_i))\}_{i=1}^{m}$.
\end{remark}

By computing~\eqref{eq:ace_cost}, ACE addresses the minimal-change
retuning problem associated with the performance constraint in~\eqref{eq:optim_problem},
but it does not enforce closed-loop stability. This motivates the uncertainty-aware stability screen developed in Section~\ref{sec:approach}, which is built on the Bayesian Sparse Identification of Nonlinear Dynamics (SINDy)~\cite{brunton_16} framework introduced next.

\subsection{SINDy}

Using closed-loop experiments (c.f.~\eqref{eq:experiments}), SINDy aims to learn the unknown map $\ve f$ as a sparse linear combination of candidate functions.

To this end, we define a library of nonlinear features
\begin{equation}\label{eq:library}
    \ve\Theta(\ve x_t, \ve u_t)
    =
    [\varphi_1(\ve x_t, \ve u_t), 
    \varphi_2(\ve x_t, \ve u_t),
    \dots,
    \varphi_{l}(\ve x_t, \ve u_t)],
\end{equation}
where each $\varphi_q$ is a candidate function of the state and input, satisfying $\varphi_q\colon\mathbb{R}^{n_x}\times\mathbb{R}^{n_u}\to\mathbb{R}$, for $q=1,\ldots,l$. Therefore, for any pair $(\ve x_t,\ve u_t)$, the evaluated library $\ve\Theta(\ve x_t,\ve u_t)$ is a row vector in $\mathbb{R}^{1\times l}$.

SINDy postulates that the unknown transition map $\ve f$ can be approximated as a linear combination of the $l$ candidate functions in~\eqref{eq:library}. Following the standard row-wise regression convention, the one-step map is written as
\begin{equation}
\label{fsindy}
    \ve f(\ve x_t,\ve u_t)^\top 
    \approx 
    \ve\Theta(\ve x_t,\ve u_t)\ve\Xi,
\end{equation}
where $\ve\Xi=[\ve\xi_1,\ldots,\ve\xi_{n_x}]\in\mathbb{R}^{l\times n_x}$ is the SINDy coefficient matrix, and $\ve\xi_j\in\mathbb{R}^{l}$ denotes the column of coefficients associated with the $j$th state equation.

Equivalently, for $\ve x_t=[x_{1,t},x_{2,t},\ldots,x_{n,t}]^\top$ and residual vector $\ve \varepsilon_t$ with components $\varepsilon_{j,t}$, the SINDy approximation yields the following component-wise one-step model:
\begin{equation}
    x_{j,t+1}
    =
    \ve\Theta(\ve x_t,\ve u_t)\ve\xi_j
    +
    \varepsilon_{j,t},
    \qquad
    j=1,\ldots,n_x,
    \label{eq:sindy_column}
\end{equation}
where the residual vector $\ve \varepsilon_t$ is assumed to be zero-mean Gaussian distributed with covariance $\ve \Sigma_w$. Consistently with the setup in~\eqref{eq:experiments}, the proposed safe tuning framework uses the $m$ closed-loop trajectories to construct one-step transition data for SINDy. For each experiment $i$ of length $N_i$, define
\begin{equation*}
    \ve X_i=
    \begin{bmatrix}
        \ve x_{0}^{\langle i\rangle\top}\\
        \ve x_{1}^{\langle i\rangle\top}\\
        \vdots\\
        \ve x_{N_i-1}^{\langle i\rangle\top}
    \end{bmatrix}
    \in\mathbb{R}^{N_i\times n_x},
    \qquad
    \ve X_i^+=
    \begin{bmatrix}
        \ve x_{1}^{\langle i\rangle\top}\\
        \ve x_{2}^{\langle i\rangle\top}\\
        \vdots\\
        \ve x_{N_i}^{\langle i\rangle\top}
    \end{bmatrix}
    \in\mathbb{R}^{N_i\times n_x}.
\end{equation*}
Here, $\ve X_i$ and $\ve X_i^+$ contain the current and one-step-ahead states. The associated library matrix is
\begin{equation*}
    \ve\Phi_{\Theta}^{\langle i\rangle}
    =
    \begin{bmatrix}
        \ve\Theta(\ve x_{0}^{\langle i\rangle},\ve u_{0}^{\langle i\rangle})\\
        \ve\Theta(\ve x_{1}^{\langle i\rangle},\ve u_{1}^{\langle i\rangle})\\
        \vdots\\
        \ve\Theta(\ve x_{N_i-1}^{\langle i\rangle},\ve u_{N_i-1}^{\langle i\rangle})
    \end{bmatrix}
    \in\mathbb{R}^{N_i\times l}.
\end{equation*}
Stacking all experiments gives
\begin{equation*}
    \ve X^+
    =
    \begin{bmatrix}
        \ve X_1^+\\
        \vdots\\
        \ve X_m^+
    \end{bmatrix}
    \in\mathbb{R}^{N\times n_x},
    \qquad
    \ve\Phi_\Theta
    =
    \begin{bmatrix}
        \ve\Phi_{\Theta}^{\langle 1\rangle}\\
        \vdots\\
        \ve\Phi_{\Theta}^{\langle m \rangle}
    \end{bmatrix}
    \in\mathbb{R}^{N \times l},
\end{equation*}
where $N=\sum_{i=1}^m N_i$. Hence, the stacked SINDy regression problem is
\begin{equation}
    \ve X^+
    =
    \ve\Phi_\Theta \ve\Xi
    +
    \ve W,
    \label{eq:sindy_stacked_regression}
\end{equation}
with $\ve W\in\mathbb{R}^{N\times n_x}$ the matrix formed by the residuals.

The key assumption in SINDy is that only a small number of library terms are needed to describe the local dynamics. This sparsity is typically enforced through sparse-regression methods such as STLSQ~\cite{brunton_16}, LASSO~\cite{tibshirani_96}, or relaxed variants~\cite{zheng_18}, helping avoid overparameterization when only limited closed-loop data are available. In the notation used below, $i$ indexes
closed-loop experiments, $t$ time samples within each experiment, $q$ library
terms, and $j$ state equations. The next section places~\eqref{eq:sindy_stacked_regression} in a Bayesian setting, so that the uncertainty in the local dynamics can be quantified and used for stability screening.
\subsection{Bayesian formulation}\label{sec:bayesian_sindy}
The SINDy regression in the previous section provides a deterministic one-step relation between the evaluated library and the state transition data. However, in the proposed safe tuning framework, the available transition data are noisy and limited, making a point estimate of the SINDy coefficients insufficient to assess how model uncertainty affects the stability of a candidate controller.

Motivated by this limitation, we introduce a Bayesian inference formulation for
the local SINDy model, following Uncertainty Quantification SINDy
(UQ-SINDy)~\cite{hirsh_22}. The goal is to infer a distribution over plausible coefficient matrices and residual uncertainty, rather than selecting a single deterministic model.

From the stacked SINDy regression model in~\eqref{eq:sindy_stacked_regression}, 
the residual matrix $\ve W$ collects the prediction errors 
$\ve \varepsilon_{t}^{\langle i\rangle}$ across all $m$ experiments, each with covariance 
$\ve\Sigma_w$. We treat the unknown coefficient matrix 
$\ve\Xi$ and the error covariance $\ve\Sigma_w$ as random variables. Under the 
Gaussian noise model defined in Section~\ref{sec:problem} and the transition model in \eqref{eq:sindy_column}, we have
\begin{equation}
    \ve x_{t+1}^{\langle i\rangle} \hspace{-0.03cm}
    \mid \hspace{-0.03cm}
    \ve x_{t}^{\langle i\rangle}\hspace{-0.02cm},\hspace{-0.02cm} \ve u_{t}^{\langle i\rangle}\hspace{-0.02cm},\ve\Xi,\ve\Sigma_w \hspace{-0.04cm}\sim \hspace{-0.04cm}
    \mathcal N 
    \hspace{-0.03cm}\left(\hspace{-0.02cm}
        (\ve\Theta(\ve x_{t}^{\langle i\rangle},\hspace{-0.03cm}\ve u_{t}^{\langle i\rangle})\ve\Xi)^{\hspace{-0.05cm}\top}\hspace{-0.02cm},
        \ve\Sigma_w \hspace{-0.02cm}
    \right)\hspace{-0.05cm}.
    \label{eq:bayesian_one_step_likelihood}
\end{equation}
Following the notation of Remark~\ref{rem:dataset}, let 
$\mathcal D$ denote the transition dataset  $\{(\ve r_{t}^{\langle i\rangle},\ve u_{t}^{\langle i\rangle},\ve x_{t}^{\langle i\rangle})\}_{t=0}^{N_i}$ collected over  $i=1,\ldots,m$ experiments. In the Bayesian setting, the objective is to infer the posterior distribution $p(\ve\Xi,\ve\Sigma_w,\ve\zeta \mid \mathcal D)$, where $\ve\zeta$ denotes auxiliary variables associated with the sparsity structure of the SINDy model~\cite[Sec.~2.3]{hirsh_22}. These variables may encode, for instance, whether each library term is active in each state equation.

Applying Bayes' rule gives
\begin{equation}
    p(\ve\Xi,\ve\Sigma_w,\ve\zeta \mid \mathcal D)
    \propto
    p(\mathcal D\mid \ve\Xi,\ve\Sigma_w)
    p(\ve\Xi\mid\ve\zeta,\ve\Sigma_w)
    p(\ve\zeta)
    p(\ve\Sigma_w).
    \label{eq:bayes_rule_sindy}
\end{equation}
The likelihood term 
$p(\mathcal D\hspace{-0.05cm}\mid\hspace{-0.05cm} \ve\Xi,\ve\Sigma_w)$ 
measures how well a given SINDy model explains the observed one-step transitions. 
The conditional prior 
$p(\ve\Xi\mid\ve\zeta,\ve\Sigma_w)$ 
encodes sparsity assumptions on the SINDy coefficients, while 
$p(\ve\zeta)$ and $p(\ve\Sigma_w)$ specify prior assumptions on the auxiliary sparsity variables and the residual covariance, respectively. The resulting posterior distribution over plausible local dynamical models is then used for stability screening in the proposed safe tuning framework described in Section~\ref{sec:approach}.

\section{Proposed Approach: UQ-SINDy Stability Screen for ACE}
\label{sec:approach}
ACE has been successfully applied both to the explanation of complex machine-learning models and to the sample-efficient retuning of linear controllers. In the latter setting, ACE seeks minimal parameter changes that turn a failing baseline controller into a passing one with few closed-loop queries. However, ACE by itself does not explicitly quantify how uncertainty in the local closed-loop dynamics affects the stability of the candidate counterfactual controllers (second constraint in~\eqref{eq:optim_problem}).

In this section, we introduce UQ-ACE, a stability-aware extension of the ACE framework. UQ-ACE augments ACE with a data-driven stability screen based on the
UQ-SINDy posterior introduced in Section~\ref{sec:bayesian_sindy}. This
posterior is used as a Bayesian local modeling tool: the closed-loop trajectories collected during ACE are used to infer a posterior distribution over plant models. Each posterior sample induces a possible closed-loop model for a candidate controller, and the candidate is accepted for testing only if it stabilizes a prescribed fraction of these posterior models. This preserves the minimal-change and sample-efficient nature of ACE while adding an uncertainty-aware stability requirement before deploying a new controller.

\subsection{SINDy posterior distribution}\label{sec:sindy_post}
We now describe the Bayesian SINDy posterior used as the uncertainty model in
UQ-ACE. The construction combines a closed-loop trajectory likelihood with a
sparsity-promoting prior on the SINDy coefficients.
\subsection*{Likelihood from closed-loop trajectories}
Starting from the Bayesian one-step model in~\eqref{eq:bayesian_one_step_likelihood}, we assume that the transition errors are conditionally Gaussian and independent across time and experiments, given the SINDy coefficients $\ve\Xi$, the residual covariance $\ve\Sigma_w$, and the closed-loop signals. Because the data are collected in closed loop, the input $\ve u_{t}^{\langle i\rangle}$ is not an externally prescribed excitation signal, but is determined by the controller and the measured state--reference history. 

Under these assumptions, the likelihood of the $M$ closed-loop trajectories is given by
\begin{align}
    p(\mathcal{D} \hspace{-0.07cm}\mid\hspace{-0.07cm} \ve\Xi,\hspace{-0.02cm}\ve\Sigma_w)
    \hspace{-0.01cm}&\!=\!\hspace{-0.02cm}
    \prod_{i=1}^{m} \!
    \hspace{-0.02cm}\prod_{t=0}^{N_i\hspace{-0.02cm}-\hspace{-0.02cm}1}\!
    \hspace{-0.04cm}\mathcal N\hspace{-0.02cm}\!\left(\!
        \ve x_{t+1}^{\langle i\rangle}; \hspace{-0.03cm}
        (
        \ve\Theta(\ve x_{t}^{\langle i\rangle}\hspace{-0.04cm},\hspace{-0.02cm}\ve u_{t}^{\langle i\rangle})\ve\Xi
        )^{\hspace{-0.04cm}\top}\hspace{-0.05cm},\hspace{-0.03cm}
        \ve\Sigma_w
    \right). 
    \label{eq:closed_loop_sindy_likelihood}
\end{align}
\subsection*{Sparsity priors and posterior sampling}
While deterministic SINDy typically enforces sparsity through sparse regression penalties or thresholding rules, the Bayesian formulation encodes the same structural assumption through a sparsity-promoting prior on $\ve\Xi$. Using the
auxiliary variables $\ve\zeta$ introduced in Section~\ref{sec:bayesian_sindy}, we write this prior as
\begin{equation}
    p(\ve\Xi,\ve\zeta,\ve\Sigma_w)
    =
    p(\ve\Xi\mid \ve\zeta,\ve\Sigma_w)
    p(\ve\zeta)
    p(\ve\Sigma_w).
    \label{eq:sindy_prior_factorization}
\end{equation}
Different choices of $\ve\zeta$ lead to different sparsity-inducing priors. For example, under a spike-and-slab prior~\cite{mitchell_88}, $\ve\zeta=\ve\lambda$, where $\lambda_{qj}\in\{0,1\}$ is a binary variable indicating whether the $q$th library term is active in the $j$th state equation. For each coefficient $\xi_{qj}$, we assign the hierarchical spike-and-slab prior
\begin{equation*}
\begin{aligned}
    \xi_{qj}\mid \lambda_{qj}
    &\sim
    (1-\lambda_{qj})\mathcal N(0,\tau_{\mathrm{spike}}^2)
    +
    \lambda_{qj}\mathcal N(0,\tau_{\mathrm{slab}}^2), \\
    \lambda_{qj} &\sim \mathrm{Bernoulli}(\pi).
\end{aligned}
\end{equation*}
with $\tau_{\mathrm{spike}}^2 \ll \tau_{\mathrm{slab}}^2$. Thus, inactive terms are strongly shrunk toward zero, while active terms are allowed to take larger values.

Alternatively, in a regularized horseshoe prior~\cite{juho_17}, the auxiliary variables are continuous shrinkage parameters rather than binary indicators. For each coefficient $\xi_{qj}$, one typical hierarchical form is
\begin{equation*}
\begin{aligned}
    \xi_{qj}\mid \lambda_{qj},\tau_j,c_j,\sigma_j^2
    &\sim
    \mathcal N\!\left(0,\sigma_j^2\tau_j^2\tilde{\lambda}_{qj}^{\,2}\right),
    \quad
    \tilde{\lambda}_{qj}^{\,2}
    =
    \frac{c_j^2\lambda_{qj}^2}{c_j^2+\tau_j^2\lambda_{qj}^2}, \\
    \lambda_{qj}
    &\sim
    \mathcal C^{+}(0,1),
    \qquad
    \tau_j \sim \mathcal C^{+}(0,\tau_0).
\end{aligned}
\label{eq:horseshoe_prior}
\end{equation*}
Here, $\mathcal C^+$ denotes the half-Cauchy distribution, $\lambda_{qj}$ is a local shrinkage parameter, $\tau_j$ is a global shrinkage parameter for the $j$th state equation, $\sigma_j^2$ denotes the residual variance of that equation, and $c_j$ is a slab-scale hyperparameter. In this case, $\ve\zeta=(\ve\lambda,\ve\tau,\ve c)$, so that the prior shrinks most coefficients toward zero while allowing selected ones to remain large when supported by the data.

Combining the closed-loop likelihood in~\eqref{eq:closed_loop_sindy_likelihood} 
with the sparsity prior in~\eqref{eq:sindy_prior_factorization} yields the Bayesian SINDy posterior in~\eqref{eq:bayes_rule_sindy}. 
In this work, we
implement both spike-and-slab and regularized horseshoe priors. For the spike-and-slab prior, the auxiliary variables are binary inclusion
indicators, which makes it natural to use a Gibbs sampler that alternates
between conditional updates of the SINDy coefficients, the sparsity variables,
and, when estimated, the residual covariance. The corresponding conditional distributions satisfy
\begin{equation*}
    \begin{aligned}
    p(\ve\Xi\mid \ve\lambda,\ve\Sigma_w,\mathcal D)
    &\propto
    p(\mathcal D\mid \ve\Xi,\ve\Sigma_w)
    p(\ve\Xi\mid \ve\lambda,\ve\Sigma_w), \\
    p(\ve\lambda\mid \ve\Xi,\ve\Sigma_w,\mathcal D)
    &\propto
    p(\ve\Xi\mid \ve\lambda,\ve\Sigma_w)
    p(\ve\lambda), \\
    p(\ve\Sigma_w\mid \ve\Xi,\ve\lambda,\mathcal D)
    &\propto
    p(\mathcal D\mid \ve\Xi,\ve\Sigma_w)
    p(\ve\Xi\mid \ve\lambda,\ve\Sigma_w)
    p(\ve\Sigma_w).
    \end{aligned}
\end{equation*}
For the regularized horseshoe prior, the sparsity variables are continuous
global and local shrinkage parameters rather than binary indicators. We therefore sample the posterior using the No-U-Turn Sampler (NUTS)~\cite{hoffman_14}, following
the UQ-SINDy formulation of \cite{hirsh_22}. 
The resulting posterior can be written as
\begin{equation*}
\begin{aligned}
    &p(\ve\Xi,\ve\lambda,\ve\tau,\ve c,\ve\Sigma_w\mid \mathcal D) \\
    &\quad\propto
    p(\mathcal D\mid \ve\Xi,\ve\Sigma_w)
    p(\ve\Xi\mid \ve\lambda,\ve\tau,\ve c,\ve\Sigma_w)
    p(\ve\lambda,\ve\tau,\ve c,\ve\Sigma_w).
\end{aligned}
\label{eq:rhs_posterior_sindy}
\end{equation*}
%
The resulting samples $\big\{\ve\Xi^{(s)},\ve\zeta^{(s)},\ve\Sigma_w^{(s)}\big\}_{s=1}^{S}$ induce $S$ plausible local SINDy models consistent with the collected closed-loop data.


\subsection{Posterior closed-loop models}
\label{subsec:posterior_closed_loop_models}

The posterior samples obtained in the previous subsection define an ensemble of plausible local plant models. For each posterior sample $\ve\Xi^{(s)}$, with $s=1,\ldots,S$, we define the corresponding SINDy dynamics as
\begin{equation*}
    \hat{\ve f}^{(s)}(\ve x,\ve u)
    =
    \ve\Theta(\ve x,\ve u)\ve\Xi^{(s)} .
\end{equation*}
The model $\hat{\ve f}^{(s)}$ represents one plausible nonlinear state-update map consistent with the closed-loop data and the sparsity-promoting posterior. Each posterior dynamics is linearized around a nominal operating point $(\bar{\ve x},\bar{\ve u})$. This operating point can be selected from prior process knowledge or estimated
from the collected trajectories.

For each posterior sample, we compute
\begin{equation*}
\begin{aligned}
    \ve A_p^{(s)}
    &=
    \left.
    \frac{\partial \hat{\ve f}^{(s)}(\ve x,\ve u)}
         {\partial \ve x}
    \right|_{(\bar{\ve x},\bar{\ve u})},
    \qquad
    \ve B_p^{(s)}
    =
    \left.
    \frac{\partial \hat{\ve f}^{(s)}(\ve x,\ve u)}
         {\partial \ve u}
    \right|_{(\bar{\ve x},\bar{\ve u})}.
\end{aligned}
\end{equation*}
Thus, each posterior sample induces the local plant model
\begin{equation*}
    \Delta \ve x_{t+1}
    =
    \ve A_p^{(s)}\Delta \ve x_t
    +
    \ve B_p^{(s)}\Delta \ve u_t,
    \qquad
    \Delta \ve y_t
    =
    \ve C_p\Delta \ve x_t .
\end{equation*}

We now combine each local posterior plant with a candidate controller parameterized by $\ve\theta$.
Substituting \eqref{eq:System_b}-\eqref{eq:System_d} in their incremental form into the plant-controller interconnection yields, for each posterior sample, a local closed-loop model of the form
\begin{equation*}
    \begin{bmatrix}
        \Delta \ve x_{t+1} \\
        \Delta \ve x_{c,t+1}
    \end{bmatrix}
    =
    \ve A_{\mathrm{cl}}^{(s)}(\ve\theta)
    \begin{bmatrix}
        \Delta \ve x_t \\
        \Delta \ve x_{c,t}
    \end{bmatrix}
    +
    \ve B_{\mathrm{cl}}^{(s)}(\ve\theta)
         \Delta \ve r_{t},
\end{equation*}
with
\begin{equation}
    \ve A_{\mathrm{cl}}^{(s)}(\ve\theta)
    =
    \begin{bmatrix}
        \ve A_p^{(s)}
        -
        \ve B_p^{(s)}
        \ve D_c(\ve\theta)
        \ve C_p
        &
        \ve B_p^{(s)}
        \ve C_c(\ve\theta)
        \\
        -
        \ve B_c(\ve\theta)
        \ve C_p
        &
        \ve A_c(\ve\theta)
    \end{bmatrix}.
    \label{eq:posterior_closed_loop_matrix}
\end{equation}
\subsection{Posterior stability probability}
\label{subsec:posterior_stability_probability}

For each posterior sample $s$, the augmented closed-loop matrix
$\ve A_{\mathrm{cl}}^{(s)}(\ve\theta)$ in~\eqref{eq:posterior_closed_loop_matrix}
is used to assess the local closed-loop stability of the candidate controller
$\ve\theta$. Given a stability margin $\epsilon_{\mathrm{stab}}>0$,
the $s$th posterior closed-loop model is classified as stable if
\begin{equation}
    \rho\!\left(
    \ve A_{\mathrm{cl}}^{(s)}(\ve\theta)
    \right)
    <
    1-\epsilon_{\mathrm{stab}},
    \label{eq:posterior_sample_stability_condition}
\end{equation}
where $\rho(\cdot)$ denotes the spectral radius. The margin $\epsilon_{\mathrm{stab}}$ makes the stability test conservative by
requiring the dominant closed-loop eigenvalue to lie strictly inside the unit
circle. This margin can be calibrated from historical closed-loop data, for
instance by examining the dominant closed-loop eigenvalues of previously
implemented controllers.

The posterior empirical probability that $\ve\theta$ stabilizes the closed-loop
system is then defined as
\begin{equation}
    \widehat p_{\mathrm{stab}}^{m}(\ve\theta)
    =
    \frac{1}{S}
    \sum_{s=1}^{S}
    \mathbf{1}
    \left\{
    \rho\!\left(
    \ve A_{\mathrm{cl}}^{(s)}(\ve\theta)
    \right)
    <
    1-\epsilon_{\mathrm{stab}}
    \right\},
    \label{eq:posterior_stability_probability}
\end{equation}
where $m$ denotes the number of closed-loop experiments collected so far and
$\mathbf{1}\{\cdot\}$ is the indicator function. Thus,
$\widehat p_{\mathrm{stab}}^{M}(\ve\theta)$ estimates the fraction of posterior
plant models for which the candidate controller satisfies the prescribed local
stability margin.
\subsection{Stability-constrained UQ-ACE acquisition}
\label{subsec:stability_constrained_uq_ace}

The previous subsections define a posterior ensemble of local closed-loop models
and the corresponding stability probability
$\widehat p_{\mathrm{stab}}^{m}(\ve\theta)$ (cf.~\eqref{eq:posterior_stability_probability}). We now use this probability as a
safety filter inside the ACE controller-tuning loop. 

Let $m$ denote the number
of closed-loop experiments collected so far. At iteration $m+1$, UQ-ACE proposes
a safe candidate controller by solving
\begin{equation}
    \ve\theta_{m+1}
    =
    \underset{\ve\theta\in\Theta_{\mathrm{trust}}}{\arg\max}
    \;
    \mathrm{EI}_{m}(\ve\theta)
    \quad
    \mathrm{s.t.}
    \quad
    \widehat p_{\mathrm{stab}}^{m}(\ve\theta)
    \ge
    1-\delta_{\mathrm{stab}} .
    \label{eq:safe_uq_ace_acquisition}
\end{equation}
Here, $\Theta_{\mathrm{trust}}$ denotes the admissible controller-parameter
region considered by ACE, $\mathrm{EI}_{m}$ is the expected improvement
acquisition function computed from the current GPC model, and
$1-\delta_{\mathrm{stab}}$ is the required posterior probability of local
closed-loop stability.

If a candidate controller has high expected improvement but does not satisfy the stability constraint in~\eqref{eq:safe_uq_ace_acquisition}, it is rejected
without being tested on the physical system. If the constraint is satisfied, a closed-loop experiment is performed with
$\ve\theta_{m+1}$. This experiment produces a new class label
$\ell_{m+1}$, used by ACE to update its classifier, and a new batch of
closed-loop data $ \mathcal D_{m+1}
    =
    \{\ve r_{t}^{\langle m+1\rangle},\ve u_{t}^{\langle m+1\rangle},\ve x_{t}^{\langle m+1\rangle}\}_{t=0}^{N_{m+1}},
$
is used to update the UQ-SINDy posterior.

The two datasets are therefore updated as $    \mathcal P
    =
    \mathcal P
    \cup
    \{(\ve\theta_{m+1},\ell_{m+1})\},$ and $    \mathcal D
    =
    \mathcal D
    \cup
    \{\mathcal D_{m+1}\}$.
Thus, each accepted experiment simultaneously refines the ACE classifier and the UQ-SINDy posterior. The resulting UQ-ACE procedure is summarized in Algorithm~\ref{alg:uq_ace_safe_tuning}.
\begin{algorithm}
\caption{UQ-ACE controller fine tuning}
\label{alg:uq_ace_safe_tuning}
\begin{algorithmic}[1]
\REQUIRE Baseline controller $\ve\theta_0$, historical closed-loop data
   $\mathcal D_{\mathrm{hist}}=\{\mathcal D_i\}_{i=1}^{N_{\mathrm{hist}}}$, $\mathcal P_{\mathrm{hist}}=\{\mathcal P_i\}_{i=1}^{N_{\mathrm{hist}}}$, $\delta_{\mathrm{stab}}$ 
\ENSURE Counterfactual controller $\ve\theta_{\mathrm{CFE}}$

\STATE Initialize $\mathcal P \leftarrow \mathcal P_{hist}$ and
$\mathcal D \leftarrow \mathcal D_{hist}$

\REPEAT
    \STATE Update ACE surrogate classifier using $\mathcal P$
    \STATE Update UQ-SINDy posterior using $\mathcal D$
    \STATE Get candidate controller $\ve\theta_k$ via~\eqref{eq:ei}
    \STATE Construct posterior closed-loop models 
    via \eqref{eq:posterior_closed_loop_matrix}
    \STATE Compute $\hat p^m_{\mathrm{stab}}(\ve\theta_{\mathrm{k}})$
    via \eqref{eq:posterior_stability_probability}

    \IF{$\hat p^m_{\mathrm{stab}}(\ve\theta_k)
    < 1-\delta_{\mathrm{stab}}$}
        \STATE Reject $\ve\theta_K$ without closed-loop testing
        \STATE\textbf{continue}
    \ENDIF
        \STATE Test $\ve\theta_k$ in closed loop and get $\ell_k$, $\mathcal D_k$
        \STATE 
        $\mathcal P \leftarrow \mathcal P \cup
        \{\ve\theta_k,\ell_k\}$
        \STATE 
        $\mathcal D \leftarrow \mathcal D \cup
        \{\mathcal D_k\}$
\UNTIL{ACE stopping criterion is satisfied}

\STATE \textbf{return} $\ve\theta_{\mathrm{CFE}}$ from the updated ACE surrogate
\end{algorithmic}
\end{algorithm}

\section{Simulation example}
\label{sec:simulations}
In this section, we evaluate the proposed stability-constrained UQ-ACE method on an industrially motivated axial-compressor control benchmark. The simulation study has two objectives. First, we analyze the proposed method under two UQ-SINDy representations, namely Case~1, which uses a full-state model, and Case~2, which uses an input-output embedded-state model. Second, we compare the resulting stability-constrained UQ-ACE procedure against the original ACE implementation, which does not include the posterior stability filter.

\subsection{Experimental setup}
\noindent\textit{System process.}
We consider an axial compressor-plenum-throttle system regulated by inlet
guide vanes (IGVs) \cite[p.~24]{moore1985theory}. Around a healthy
operating point, the local compressor state is defined as
$   \ve x_t
    =
        [\widetilde{\Phi}_t,
        \widetilde{\Psi}_t]^\top
    \in\mathbb R^2,
$
where $\widetilde{\Phi}_t$ and $\widetilde{\Psi}_t$ denote deviations of the
flow and pressure-rise coefficients from their operating-point values. The IGV
command is denoted by $u_t\in\mathbb R$, and the control objective is to
regulate the pressure-rise output $y_t=\widetilde{\Psi}_t$ to a reference
$r_t$.

We consider two UQ-SINDy representations for the compressor system. In Case~1,
we use a full-state representation, where the dynamics are modeled as
\begin{equation*}
    \ve x_{t+1}
    =
    \ve F_{\ve\Xi}(\ve x_t,u_t)
    +
    \ve w_t,
    \qquad
    \ve w_t
    \sim
    \mathcal N(\ve 0,\ve\Sigma_w).
\end{equation*}
In Case~2, we use an embedded-state representation based only on the controlled
output, denoted by $y_t=\widetilde{\Psi}_t$. The embedded state is defined as
\begin{equation*}
    \ve x_t
    =
        [y_t,
        y_{t-1},
        u_{t-1},
        u_{t-2}]^\top,
\end{equation*}
with dynamics
\begin{equation*}
    \ve x_{t+1}
    =
    \ve S\ve x_t
    +
    \ve G u_t
    +
    \ve e_1 F_{\xi}(\ve x_t,u_t)
    +
    \ve e_1 w_t,
    \qquad
    w_t\sim\mathcal N(0,\sigma^2),
\end{equation*}
where
\begin{equation*}
    \ve S
    =
    \begin{bmatrix}
        0 & 0 & 0 & 0\\
        1 & 0 & 0 & 0\\
        0 & 0 & 0 & 0\\
        0 & 0 & 1 & 0
    \end{bmatrix},
    \qquad
    \ve G
    =
    \begin{bmatrix}
        0\\
        0\\
        1\\
        0
    \end{bmatrix},
    \qquad
    \ve e_1
    =
    \begin{bmatrix}
        1\\
        0\\
        0\\
        0
    \end{bmatrix}.
\end{equation*}

\noindent \textit{Controller.} The controller is implemented as a PID with derivative filtering, i.e.,
\[
C(s)\;=\;K_p\;+\;\frac{K_i}{s}\;+\;\frac{K_d\, s}{1+T_f s}.
\]
Accordingly, UQ-ACE searches for a four-parameter retuning, with $\ve{\theta}=[K_p, K_i, K_d, T_f]^\top$.

\noindent \textit{Parameters.} Unless otherwise noted, the compressor– plenum constants follow those in \cite{guerrero_26_ifac}. We use a normalized step of amplitude $r_{\text{step}}=1$. Measurement noises are zero-mean Gaussian with variances
$\sigma^2=0.01$, corresponding to a input signal-to-noise ratio of $20$ $[\textnormal{dB}]$. The ACE hyperparameters follow  \cite[Appendix G.1]{guerrero_26}.

\noindent \textit{Closed-loop criteria.}
We define the pass/fail function as
\begin{equation}
h(\boldsymbol{\theta}) =
\begin{cases}
1, & \text{if } |e_{\mathrm{ss}}|\le 0.01,\;
t_{\mathrm{rise}(10\to90)}\le 20\,\mathrm{s}, \\[-0.2ex]
   & \qquad t_{\mathrm{settle}}\le 80\,\mathrm{s},\;
\mathrm{OS}\le 25\%,\\
0, & \text{otherwise.}
\end{cases}
\label{eq:criteria}
\end{equation}
where $e_{\mathrm{ss}}$ is the steady-state error, $t_{\mathrm{rise}(10\to90)}$ is the rise time between the first crossings of $10\%$ and $90\%$ of the final step change, $t_{\mathrm{settle}}$ is the first time for which $|y(t)-y_\infty|\le 0.05$ for all $t\ge t_{\mathrm{settle}}$, and $\mathrm{OS}$ is the percent of overshoot.

\noindent\textit{UQ-SINDy setting.}
For both representations, the nonlinear maps are parameterized using the same
polynomial SINDy library $\ve\Theta(\ve x_t,u_t)$, including constant, linear,
quadratic, and input-coupling terms. Sampled SINDy models are linearized
around a stable compressor operating point, following the operating-point
notation of the Moore--Greitzer compressor characteristic
\cite[Fig.~3.1]{moore1985theory}. In the main experiments, the posterior over
the SINDy coefficients is computed using the spike-and-slab prior with inclusion probability
$\pi_{\mathrm{incl}}=0.25$ and a Gaussian residual noise model. For each
UQ-SINDy update, we draw $S=1000$ posterior samples after a burn-in period of
200 Gibbs iterations. The required posterior stability probability is set to
$1-\delta_{\mathrm{stab}}=0.97$ for both cases.

\subsection{Test 1: Fully observable}
Before comparing UQ-ACE and ACE over multiple runs, we assess the sensitivity of
the UQ-SINDy stability screen to the sparsifying prior. We consider the
representative Case~1 instance
$\ve\theta_0=[3.1283,\,0.0268,\,7.5788,\,0.0885]^\top$ and fit both
spike-and-slab and regularized horseshoe posteriors using the same closed-loop
training data, SINDy library, and ACE-generated candidate controller. For the prior-sensitivity comparison in Table~\ref{tab:prior_sensitivity_case1},
we use $S=100$ posterior samples after 50 burn-in iterations for each prior, as
the goal is only to assess the qualitative effect of the sparsifying prior.

\begin{table}[t]
    \centering
    \caption{Prior-sensitivity analysis for one representative Case~1 instance.}
    \label{tab:prior_sensitivity_case1}
    \scriptsize
    \setlength{\tabcolsep}{3pt}
    \resizebox{\columnwidth}{!}{
    \begin{tabular}{l|c|c|c|c|c|c}
        \hline
        Prior 
        & $\hat p_0$
        & $\rho_0^{\mathrm{med}}$
        & $\hat p_{\mathrm{CFE}}$
        & $\rho_{\mathrm{CFE}}^{\mathrm{med}}$
        & $W_{90}$
        & Time \\
        \hline
        Spike-slab
        & 0.00
        & 0.9993
        & 0.60
        & 0.9990
        & $7.6{\times}10^{-4}$
        & 1.8 min. \\
        Reg. horse.
        & 0.00
        & 0.9993
        & 0.23
        & 0.9990
        & $2.8{\times}10^{-3}$
        & 67.2 min. \\
        \hline
    \end{tabular}}
\end{table}
Both priors assign zero posterior stability probability to the baseline
controller under the prescribed stability margin, as shown in
Table~\ref{tab:prior_sensitivity_case1}. The median closed-loop spectral radii are nearly identical
for the two priors, both for $\boldsymbol{\theta}_0$ and
$\boldsymbol{\theta}_{\mathrm{CFE}}$, suggesting that the nominal local dynamics
identified by the two posterior models are similar. Here, $W_{90}$ denotes the
width of the 90\% posterior interval of the closed-loop spectral radius,
computed as $q_{95}-q_{05}$. For the counterfactual controller, however, the
regularized horseshoe gives a more conservative stability assessment, reducing
$\hat p_{\mathrm{stab}}(\boldsymbol{\theta}_{\mathrm{CFE}})$ from $0.60$ to
$0.23$, consistently with its larger posterior uncertainty as indicated by
$W_{90}$. In the remaining experiments, we select the spike-and-slab prior as
the default prior because it provides an explicit sparsity representation
through binary inclusion variables and has substantially lower computational time.



Table~\ref{tab:case1_full_state} summarizes results over 100 ACE runs for three fixed instances $\boldsymbol{\theta}_{01:03}$, comparing UQ-ACE and ACE initialized with 30 historical controllers. Each entry reports the mean fraction of runs that found a counterfactual (\emph{Validity}), the mean number of additional closed-loop tests (\emph{Iters.}), and the mean component-wise parameter shift $(\Delta K_p;\Delta K_i;\Delta K_d;\Delta T_f)$, where $\Delta=\boldsymbol{\theta}_{\text{CFE}}-\boldsymbol{\theta}_0$, so that positive values indicate an increase with respect to the baseline. We also report the average number of candidates rejected by the UQ-SINDy stability screen (\emph{Screened out}) and the average number of invalid closed-loop trajectories
(\emph{Invalid traj.}), corresponding to candidate controllers deployed on the simulated closed loop
that resulted in instability; values below one
indicate occasional invalid deployments across the 100 runs. Overall, the average parameter changes are slightly smaller for UQ-ACE, although both methods produce broadly similar retunings. The main difference appears in the additional closed-loop tests: the UQ-SINDy safety filter reduces the number of deployed experiments and prevents unstable candidates from being tested. In the reported cases, UQ-ACE produces no invalid trajectories, whereas ACE occasionally deploys unstable candidates.



\begin{figure}
\centering
\includegraphics[trim=0mm 0mm 0mm 0mm,clip, width=1\columnwidth]{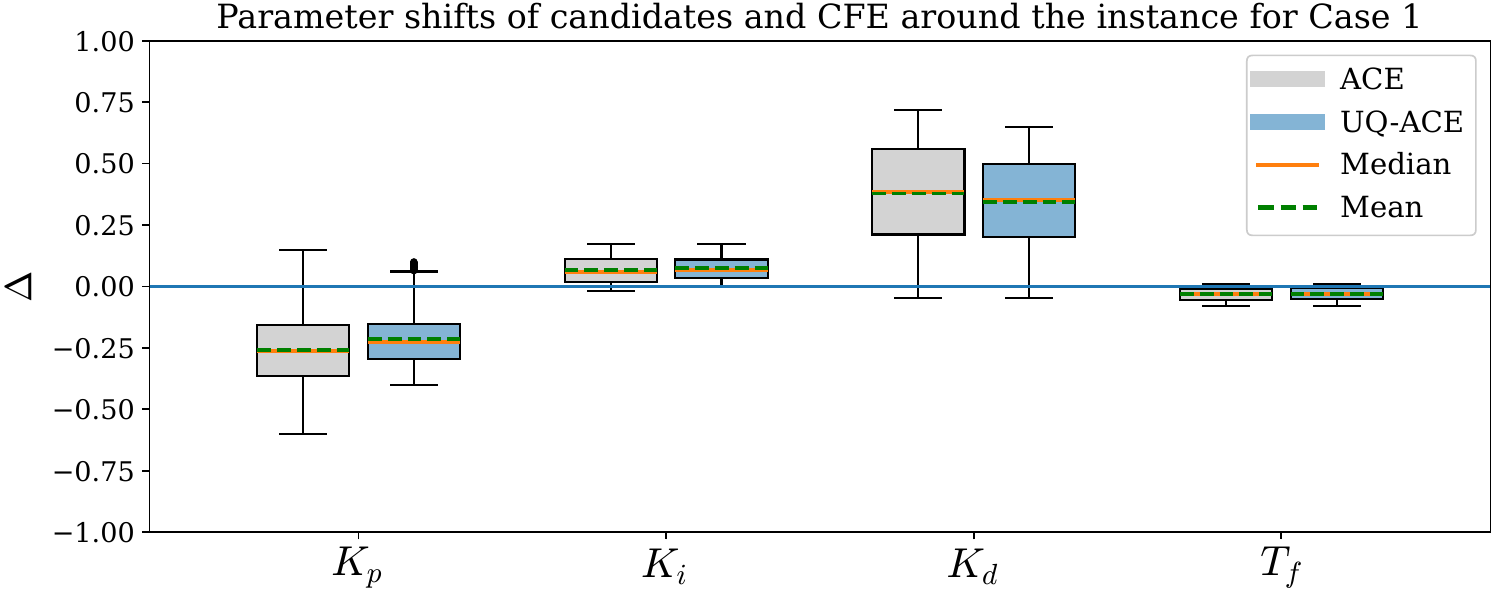}
\vspace{-0.5cm}
\caption{Case 1 - Box plot of component-wise shifts $\ve \Delta=\boldsymbol{\theta}_{\text{cand}}-\boldsymbol{\theta}_0$ across 100 ACE runs.}
    \label{fig:box_plot_single} 
\end{figure}

\begin{figure}
\centering
\includegraphics[trim=0mm 0mm 0mm 0mm,clip, width=1\columnwidth]{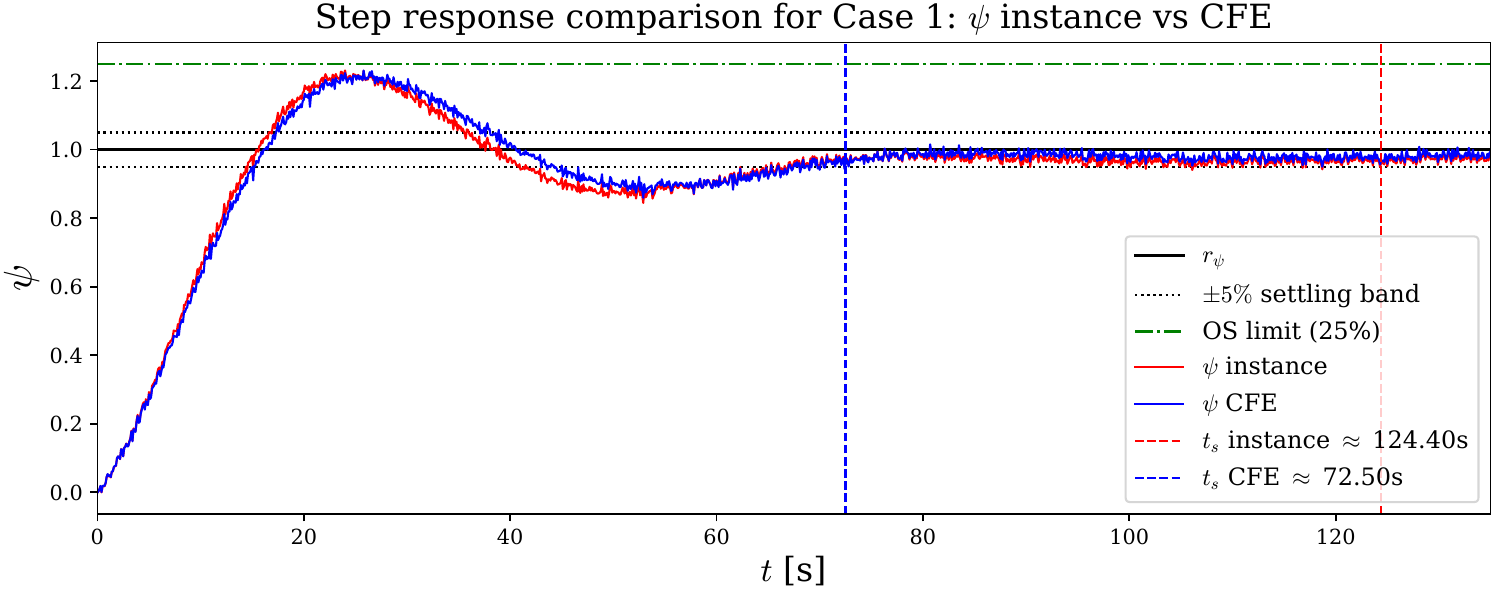}
\vspace{-0.6cm}
\caption{Case 1 - Step responses: baseline $\boldsymbol{\theta}_0$ (red) vs. one representative CFE (blue).}
    \label{fig:step_one} 
\end{figure}

\begin{table}[t]
    \centering
    \caption{Case 1, full-state representation: summary over 100 runs per instance.}
    \label{tab:case1_full_state}
    \scriptsize
    \setlength{\tabcolsep}{3pt}
    \renewcommand{\arraystretch}{1.05}
    \resizebox{\columnwidth}{!}{
    \begin{tabular}{l|l|c|c|c|c|c}
        \hline
        \makecell{Inst.} & \makecell{Method} & \makecell{Val.} & \makecell{Iters.} &
        \makecell{Screened\\out} & \makecell{Invalid\\traj.} &
        \makecell{$\Delta K_p;\Delta K_i;\Delta K_d;\Delta T_f$} \\
        \hline
        \multirow{2}{*}{$\boldsymbol{\theta}_{01}$}
        & UQ-ACE & 0.99 & 20 & 7.61 & 0 &
        $-0.2760;\ +0.0118;\ +0.6071;\ -0.0283$ \\
        & ACE    & 0.95 & 34 & --   & 0.37 &
        $-0.2845;\ +0.0117;\ +0.6734;\ -0.0329$ \\[1ex]

        \multirow{2}{*}{$\boldsymbol{\theta}_{02}$}
        & UQ-ACE & 1.00 & 21 & 8.83 & 0 &
        $-0.2721;\ -0.0003;\ +0.1244;\ -0.0224$ \\
        & ACE    & 0.95 & 30 & --   & 0.42 &
        $-0.2846;\ -0.0007;\ +0.1461;\ -0.0229$ \\[1ex]

        \multirow{2}{*}{$\boldsymbol{\theta}_{03}$}
        & UQ-ACE & 0.98 & 19 & 8.68 & 0 &
        $-0.3205;\ -0.0011;\ +0.1651;\ -0.0241$ \\
        & ACE    & 0.97 & 32 & --   & 0.22 &
        $-0.3750;\ -0.0011;\ +0.2149;\ -0.0250$ \\
        \hline
    \end{tabular}}
\end{table}
For the first instance, $\boldsymbol{\theta}_0=\ve\theta_{01}=[3.1283,0.0268,7.5788,$ $0.0885]$, Figure~\ref{fig:box_plot_single} shows the shifts $\Delta=\boldsymbol{\theta}_{\text{cand}}-\boldsymbol{\theta}_0$ pooled over all candidate evaluations across 100 runs. Unlike Table~\ref{tab:case1_full_state}, which reports only the final CFE averages, the figure shows the full search distribution; boxes denote the interquartile range, orange/green lines the median/mean, whiskers the non-outlier range, and circles the outliers. ACE and UQ-ACE follow similar retuning directions, but UQ-ACE yields slightly more concentrated shifts, suggesting that the UQ-SINDy filter removes aggressive candidates. In both cases, the main changes occur in $K_p$ and $K_d$, while $K_i$ and $T_f$ remain close to zero due to sparsity. Figure~\ref{fig:step_one} compares $\ve\theta_0$ with one representative CFE, $[2.869,0.0340,7.8497,0.0677]^\top$, corresponding to $\Delta=[-0.2592,0.0072,0.2709,-0.0207]^\top$. These changes reduce the settling time and turn the controller from fail to pass according to~\eqref{eq:criteria}.


\subsection{Test 2: Embedding}
For Case~2, we compare UQ-ACE against the same three ACE instances reported
in Table 1 of~\cite{guerrero_26_ifac}. Table~\ref{tab:case2_embedding} shows
that UQ-ACE achieves comparable validity with slightly fewer additional
closed-loop tests, since the UQ-SINDy stability filter rejects some candidate
controllers before deployment. As reported in~\cite{guerrero_26_ifac}, ACE
does not produce unstable trajectories in these instances, so the main role of
UQ-ACE is to add a robustness layer rather than substantially change the final
retuning.
\begin{table}[t]
    \centering
    \caption{Case 2, embedded-state representation: summary over 100 runs per instance.}
    \label{tab:case2_embedding}
    \scriptsize
    \setlength{\tabcolsep}{3pt}
    \renewcommand{\arraystretch}{1.05}
    \resizebox{\columnwidth}{!}{
    \begin{tabular}{l|l|c|c|c|c|c}
        \hline
        \makecell{Inst.} & \makecell{Method} & \makecell{Val.} &
        \makecell{iters.} & \makecell{Screened\\out} &
        \makecell{invalid\\traj.} &
        \makecell{$\Delta K_p;\Delta K_i;\Delta K_d;\Delta T_f$} \\
        \hline
        \multirow{2}{*}{$\boldsymbol{\theta}_{01}$}
        & UQ-ACE & 1.00 & 10 & 4.07 & 0 &
        $+0.0537;\ -0.1017;\ +0.0245;\ +0.0339$ \\
        & ACE    & 1.00 & 12 & --   & 0 &
        $+0.0334;\ -0.0986;\ +0.0251;\ +0.0331$ \\[1ex]

        \multirow{2}{*}{$\boldsymbol{\theta}_{02}$}
        & UQ-ACE & 1.00 & 12 & 2.28 & 0 &
        $+0.0049;\ +0.0475;\ +0.0031;\ +0.0332$ \\
        & ACE    & 1.00 & 12 & --   & 0 &
        $+0.0057;\ +0.0482;\ +0.0029;\ +0.0333$ \\[1ex]

        \multirow{2}{*}{$\boldsymbol{\theta}_{03}$}
        & UQ-ACE & 1.00 & 17 & 3.18 & 0 &
        $+0.0204;\ +0.0590;\ +0.0089;\ +0.0321$ \\
        & ACE    & 1.00 & 18 & --   & 0 &
        $+0.0230;\ +0.0599;\ +0.0076;\ +0.0324$ \\
        \hline
    \end{tabular}}
\end{table}



\begin{figure}
\centering
\includegraphics[trim=0mm 0mm 0mm 0mm,clip, width=1\columnwidth]{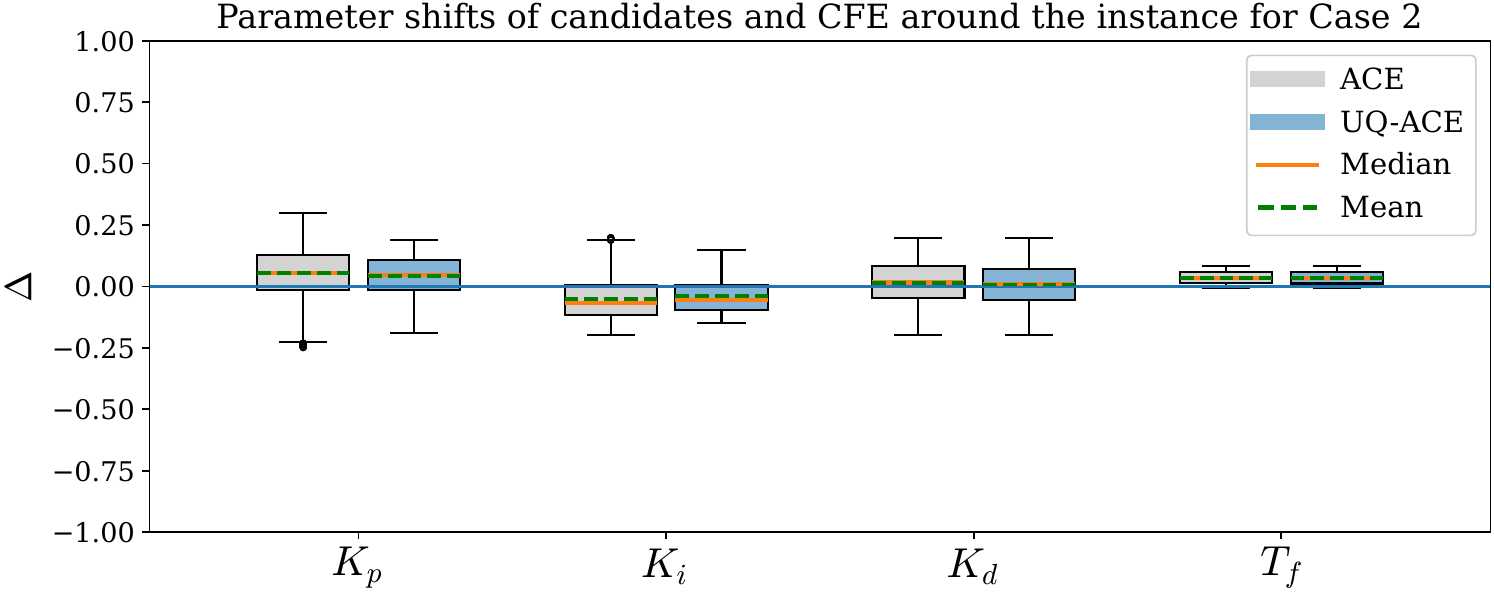}
\caption{Case 2 - Box plot of component-wise shifts $\ve \Delta=\boldsymbol{\theta}_{\text{cand}}-\boldsymbol{\theta}_0$ across 100 ACE runs.}
    \label{fig:emb_boxplot} 
\end{figure}

Figure~\ref{fig:emb_boxplot} shows similar parameter-change distributions for
ACE and UQ-ACE, with the latter slightly more concentrated around zero. This suggests
that ACE is already effective in this embedded-state setting, while UQ-ACE
preserves similar performance and screens potentially destabilizing candidates
before they are tested.
\section{Conclusions}
\label{sec:conclusions}
We have presented UQ-ACE, a stability-aware extension of ACE for
counterfactual controller tuning. The method combines ACE's sample-efficient
search with a Bayesian SINDy stability screen, using closed-loop data to
quantify uncertainty in local dynamical models and reject candidates whose
posterior stability probability falls below a prescribed threshold. Simulation
results on an axial-compressor--plenum--throttle system show that UQ-ACE
achieves valid counterfactual retunings with a number of successful experiments
comparable to ACE, while preventing unstable candidates from being deployed.
Future work includes improved unstable-loop detection and convergence-based
diagnostics, as well as extensions beyond PID tuning to nonlinear controllers.

\bibliography{ifacconf}

@article{hirsh_22,
  title={Sparsifying priors for {Bayesian} uncertainty quantification in model discovery},
  author={Hirsh, Seth M and Barajas-Solano, David A and Kutz, J Nathan},
  journal={Royal Society Open Science},
  volume={9},
  number={2},
  pages={211823},
  year={2022},
}

@InProceedings{guerrero_26,
  title = 	 {{ACE: Adapting sampling for Counterfactual Explanations}},
  author =       {Guerrero, Margarita A. and Rojas, Cristian R.},
  booktitle = 	 {Proceedings of The 8th Annual Learning for Dynamics and Control Conference},
  pages = 	 {242--264},
  year = 	 {2026},
  volume = 	 {331},
}

@article{guerrero_26_ifac,
title = {Sample-efficient counterfactual tuning for compressor pressure control},
journal = {Journal of Process Control},
volume = {162},
pages = {103723},
year = {2026},
author = {Margarita A. Guerrero and Rodrigo A. González and Cristian R. Rojas}
}

@article{hoffman_14,
  title={The {No-U-Turn} sampler: adaptively setting path lengths in {Hamiltonian Monte Carlo}},
  author={Hoffman, Matthew D and Gelman, Andrew},
  journal={J. Mach. Learn. Res.},
  volume={15},
  number={1},
  pages={1593--1623},
  year={2014}
}

@incollection{crompton2021,
  title={Data Management from the {DCS} to the {Historian}},
  author={Crompton, J.},
  booktitle={Machine Learning and Data Science in the Oil and Gas Industry},
  editor={Patrick Bangert},
  pages={83--110},
  year={2021}
}

@article{mitchell_88,
  title={Bayesian variable selection in linear regression},
  author={Mitchell, Toby J and Beauchamp, John J},
  journal={Journal of the american statistical association},
  volume={83},
  number={404},
  pages={1023--1032},
  year={1988}
}

@article{juho_17,
author = {Juho Piironen and Aki Vehtari},
title = {{Sparsity information and regularization in the horseshoe and other shrinkage priors}},
volume = {11},
journal = {Electronic Journal of Statistics},
number = {2},
pages = {5018--5051},
year = {2017},
}

@article{zheng_18,
  title={A unified framework for sparse relaxed regularized regression: {SR3}},
  author={Zheng, Peng and Askham, Travis and Brunton, Steven L and Kutz, J Nathan and Aravkin, Aleksandr Y},
  journal={IEEE Access},
  volume={7},
  pages={1404--1423},
  year={2018},
}

@article{tibshirani_96,
  title={Regression shrinkage and selection via the lasso},
  author={Tibshirani, Robert},
  journal={Journal of the Royal Statistical Society Series B: Statistical Methodology},
  volume={58},
  number={1},
  pages={267--288},
  year={1996},
  publisher={Oxford University Press}
}

@article{brunton_16,
  title={Discovering governing equations from data by sparse identification of nonlinear dynamical systems},
  author={Brunton, Steven L and Proctor, Joshua L and Kutz, J Nathan},
  journal={Proceedings of the national academy of sciences},
  volume={113},
  number={15},
  pages={3932--3937},
  year={2016}
}

@book{aastrom_05,
  title={Advanced PID control},
  author={Astr{\"o}m, K. J. and H{\"a}gglund, T.},
  year={2005},
  publisher={John Wiley \& Sons}
}

@article{Wabersich_23,
  title={Data-driven safety filters: Hamilton-jacobi reachability, control barrier functions, and predictive methods for uncertain systems},
  author={Wabersich, Kim P and Taylor, Andrew J and Choi, Jason J and Sreenath, Koushil and Tomlin, Claire J and Ames, Aaron D and Zeilinger, Melanie N},
  journal={IEEE Control Systems Magazine},
  volume={43},
  number={5},
  pages={137--177},
  year={2023},
}

@inproceedings{berkenkamp2016safe,
  title     = {Safe Controller Optimization for Quadrotors with {Gaussian} Processes},
  author    = {Berkenkamp, Felix and Schoellig, Angela P. and Krause, Andreas},
  booktitle = {Proceedings of the IEEE International Conference on Robotics and Automation},
  pages     = {493--496},
  year      = {2016},
}

@article{hewing_20,
  title={Learning-based model predictive control: Toward safe learning in control},
  author={Hewing, Lukas and Wabersich, Kim P and Menner, Marcel and Zeilinger, Melanie N},
  journal={Annual Review of Control, Robotics, and Autonomous Systems},
  volume={3},
  number={1},
  pages={269--296},
  year={2020},
}

@inproceedings{sui2015safe,
  title     = {Safe Exploration for Optimization with Gaussian Processes},
  author    = {Sui, Yanan and Gotovos, Alkis and Burdick, Joel W. and Krause, Andreas},
  booktitle = {Proceedings of the 32nd International Conference on Machine Learning},
  pages     = {997--1005},
  year      = {2015},
}

@article{khosravi_23,
  title={Safety-aware cascade controller tuning using constrained {Bayesian} optimization},
  author={Khosravi, Mohammad and K{\"o}nig, Christopher and Maier, Markus and Smith, Roy S. and Lygeros, John and Rupenyan, Alisa},
  journal={IEEE Transactions on Industrial Electronics},
  volume={70},
  number={2},
  pages={2128--2138},
  year={2023},
}

@techreport{moore1985theory,
  title={A theory of post-stall transients in multistage axial compression systems},
  author={Moore, F. K. and Greitzer, E. M.},
  year={1985},
  institution={NASA}
}

@article{wachter2018counterfactual,
  author  = {Wachter, Sandra and Mittelstadt, Brent and Russell, Chris},
  title   = {Counterfactual explanations without opening the black box: Automated decisions and the {GDPR}},
  journal = {Harvard Journal of Law \& Technology},
  volume  = {31},
  number  = {2},
  pages   = {841--887},
  year    = {2018}
}

@book{rasmussen2006gaussian,
  author    = {Rasmussen, Carl Edward and Williams, Christopher K. I.},
  title     = {Gaussian Processes for Machine Learning},
  publisher = {MIT Press},
  year      = {2006}
}

@book{Ljung1999SystemUser,
    title = {{System Identification: Theory For The User}, 2nd Ed},
    year = {1999},
    author = {Ljung, L.},
    publisher = {Prentice-Hall}
}

@book{luenberger2008linear,
  author    = {Luenberger, David G. and Ye, Yinyu},
  title     = {Linear and Nonlinear Programming, 3rd Ed},
  publisher = {Springer},
  address   = {New York},
  year      = {2008}
}

@inproceedings{breunig2000lof,
  author    = {Breunig, Markus M. and Kriegel, Hans-Peter and Ng, Raymond T. and Sander, J{\"o}rg},
  title     = {{LOF}: Identifying density-based local outliers},
  booktitle = {Proceedings of the 2000 ACM SIGMOD International Conference on Management of Data},
  pages     = {93--104},
  year      = {2000}
}

@article{jones1998efficient,
  author  = {Jones, Donald R. and Schonlau, Matthias and Welch, William J.},
  title   = {Efficient global optimization of expensive black-box functions},
  journal = {Journal of Global Optimization},
  volume  = {13},
  number  = {4},
  pages   = {455--492},
  year    = {1998}
}

@article{byrd1995limited,
  author  = {Byrd, Richard H. and Lu, Peihuang and Nocedal, Jorge and Zhu, Ciyou},
  title   = {A limited memory algorithm for bound constrained optimization},
  journal = {SIAM Journal on Scientific Computing},
  volume  = {16},
  number  = {5},
  pages   = {1190--1208},
  year    = {1995}
}

@article{baumann2025safety,
  title={Safety and optimality in learning-based control at low computational cost},
  author={Baumann, Dominik and Kowalczyk, Krzysztof and Rojas, Cristian R and Tiels, Koen and Wachel, Pawe{\l}},
  journal={IEEE Transactions on Automatic Control},
  volume={70},
  number={11},
  pages={7275--7286},
  year={2025},
}
\vspace{-2em}



\end{document}